\documentclass[a4wide]{article}

\usepackage[dvipsnames]{xcolor}
\usepackage{url}
\usepackage{algorithm}
\usepackage{algpseudocode}
\usepackage{todonotes}
\usepackage{xstring}

\usepackage{pdfpages}
\usepackage[a4paper]{geometry}
\usepackage{natbib}
\graphicspath{{figures/}}

\definecolor{darkpastelgreen}{rgb}{0.01, 0.75, 0.24}
\definecolor{bananayellow}{rgb}{1.0, 0.88, 0.21}

\def\RR{\textsf{R}\/}

\let\pkg=\strong
\newcommand{\mtc}{\mathcal}
\def\cC{{\mtc{C}}}
\newcommand{\set}[1]{\left\{#1\right\}}

\newcommand{\abs}[2][a]{
\IfEqCase{#1}{%
{a}{\left\vert#2\right\rvert}%
{0}{\vert#2\rvert}%
{1}{\big\vert#2\big\rvert}%
{2}{\Big\vert#2\Big\rvert}%
{3}{\bigg\vert#2\bigg\rvert}%
{4}{\Bigg\vert#2\Bigg\rvert}%
}[\PackageError{inner}{Undefined option to abs: #1}{}]%
}

\usepackage{natbib}
\usepackage{amsmath}
\usepackage{amssymb}
\usepackage[colorlinks=true]{hyperref}

\begin{document}

\title{Adjacency-constrained hierarchical clustering of a band similarity matrix
with application to genomics}

\author{Christophe Ambroise$^{1}$, Alia Dehman$^{2}$, Pierre Neuvial$^{3}$, \\
Guillem Rigaill$^{1,4}$ and Nathalie Vialaneix$^{5}$}

\date{\footnotesize
$^1$ Laboratoire de Mathématiques et Modélisation d'Evry, UMR CNRS 8071, Université d'Evry Val d'Essonne, 23 boulevard de France, 91037 Evry, France.\\
$^2$ Hyphen-stat, 195 Route d'Espagne, 31036 Toulouse, France.\\
$^3$ Institut de Mathématiques de Toulouse, UMR5219 CNRS, Université de Toulouse, UPS IMT, F-31062 Toulouse Cedex 9, France.\\
$^4$ Institute of Plant Sciences Paris Saclay IPS2, CNRS, INRA, Gif sur Yvette, France.\\
$^5$ MIAT, Université de Toulouse, INRA, Castanet-Tolosan, France.\\
}

\maketitle

  \begin{abstract}
    \textbf{Motivation:} Genomic data analyses such as Genome-Wide Association
    Studies (GWAS) or Hi-C studies are often faced with the problem of
    partitioning chromosomes into successive regions based on a similarity
    matrix
    of high-resolution, locus-level measurements. An intuitive way of doing this is to perform a modified Hierarchical Agglomerative Clustering (HAC), where only adjacent clusters (according to the ordering of positions within a chromosome) are allowed to be merged. But a major practical drawback of this method is its quadratic time and space complexity in the number of loci, which is  typically of the order of $10^4$ to $10^5$ for each chromosome.\\
    \textbf{Results:} By assuming that the similarity between physically distant
    objects is negligible, we are able to propose an implementation of this
    adjacency-constrained HAC with quasi-linear complexity. Our illustrations on
    GWAS and Hi-C datasets demonstrate the relevance of this assumption, and
    show that this method highlights biologically meaningful signals. Thanks to
    its small time and memory footprint, the method can be run on a standard
    laptop
    in minutes or even seconds.\\
    \textbf{Availability and Implementation}: Software and sample data are
    available as an \RR{} package, \pkg{adjclust}, that can be downloaded from
    the Comprehensive R
    Archive Network (CRAN).\\
    \textbf{Contact:} pierre.neuvial@math.univ-toulouse.fr 
\end{abstract}

\section{Introduction}

Genetic information is coded in long strings of DNA organised in
chromosomes. High-throughput sequencing such as   RNAseq, DNAseq, ChipSeq and
Hi-C makes it possible to study biological phenomena along the entire genome at a very
high resolution \citep{reuter_etal_MC2015}.

In most cases, we expect neighboring positions to be statistically dependent.
Using this \textit{a priori} information is one way of addressing the complexity of
genome-wide analyses. For instance, it is common practice to partition each
chromosome into regions, because such regions hopefully correspond to biological
relevant or interpretable units (such as genes or binding sites) and because
statistical modelling and inference are simplified at the scale of an individual region. In simple
cases, such regions are given (for example, in RNAseq analysis, only
genic and intergenic regions are usually considered and differential analysis is commonly
performed at the gene or transcript level). However, in more complex cases,
regions of interest are unknown and need to be discovered by mining the data. This
is the case in the two leading examples considered in this paper. In the context of Genome Wide Association Studies (GWAS),
region-scale approaches taking haplotype blocks into account can result in
substantial statistical gains \citep{gabriel_etal_S2002}. Hi-C studies
\citep{dixon_etal_N2012} have demonstrated the existence of topological domains,
which are megabase-sized local chromatin interaction domains correlating with
regions of the genome that constrain the spread of heterochromatin.
Hence, the problem of partitioning a chromosome into biologically
relevant regions based on measures of similarity between pairs of individual
loci has been extensively studied for genomic applications.

Recovering the ``best'' partition of $p$ loci for each possible number, $K$, of
classes is equivalent to a segmentation problem (also known as ``multiple
changepoint problem''). In the simplest scenario where the signals to be segmented
are piecewise-constant, such as in the case of DNA copy numbers in cancer
studies, segmentation can be cast as a least squares minimization problem
\citep{picard_etal_BMCB2005,hocking_etal_BMCB2013}. More generally, kernel-based
segmentation methods have been developed to perform segmentation on data
described by a similarity measure
\citep{harchaoui_cape_WSSP2007,arlot_etal_p2016}. Such segmentation problems are
combinatorial in nature, as the number of possible segmentations of $p$ loci
into $K$ blocks (for a given $K=1 \dots p$) is
${p \choose K} = \mathcal{O}(p^K)$. The ``best'' segmentation for all
$K=1 \dots p$ can be recovered efficiently in a quadratic time and space
complexity using dynamic programming. As discussed in \cite{celisse_etal_p2007},
in the case of kernel-based segmentation, this complexity cannot be improved
without making additional assumptions on the kernel (or the corresponding
similarity). Indeed, for a generic kernel, even computing the loss (that is, the
least square error) of any given segmentation in a fixed number of segments $K$
has a computational cost of $\mathcal{O}(p^2)$.

The goal of this paper is to develop heuristics that can be applied to genomic
studies in which the number of loci is so large (typically of the order of
$p=10^4$ to $10^6$) that algorithms of quadratic time and space complexity
cannot be applied. This paper stems from a modification of the classical
hierarchical agglomerative clustering (HAC)
\citep{kaufman_rousseeuw_FGDICA2009}, where only adjacent clusters (in the sense
of a pre-defined ordering) can be merged. This simple constraint is well suited
to genomic applications, in which loci can be ordered along chromosomes provided that an assembled genome is available. The resulting method can be seen as a
heuristic for segmentation; it provides not only a single partition of
the original loci, but a sequence of nested partitions.

This idea of incorporating such constraints was previously mentioned by
\citet{lebart_CAD1978} and \citet{grimm_CG1987}, and an \RR{} package
implementing this algorithm, \pkg{rioja} \citep{rioja2015}, has been
developed\footnote{available on CRAN at
  \url{https://cran.r-project.org/package=rioja}, but the package is currently
  orphaned.}. However, the algorithm remains quadratic in both time and
space. Its time complexity cannot be improved because all of the $p^2$
similarities have to be computed to perform the clustering. To circumvent this
difficulty, we assume that the similarity between physically distant loci is
zero, where two loci are deemed to be ``physically distant'' if they are separated by
more than $h$ other loci. The main contribution of this paper is to propose
an adjacency-constrained clustering algorithm with quasi-linear complexity
(namely, $\mathcal{O}(ph)$ in space and $\mathcal{O}(p(h + \log(p)))$ in time)
under this assumption, and to demonstrate its relevance for genomic studies.

The rest of the paper is organized as follows. In Section \ref{sec:method} we
describe the algorithm, its time and space complexity and its
implementation.  The resulting segmentation method is then applied on GWAS
datasets (Section~\ref{sec:appli-gwas}) and on Hi-C datasets
(Section~\ref{sec:appli-hi-c}), in order to illustrate that the above assumption
makes sense in such studies, and that the proposed methods can be used to recover
biologically relevant signals.

\section{Method}
\label{sec:method}

\subsection{Adjacency-constrained HAC with Ward's linkage}
\label{sec:adj-constr-hac}

In its unconstrained version, HAC starts with a trivial clustering where each
object is in its own cluster and iteratively merges the two most similar
clusters according to a distance function $\delta$ called a linkage criterion.
We focus on Ward's linkage, which was defined for clustering objects $(x_i)_i$ taking
values in the Euclidean space $\mathbb{R}^d$. Formally, Ward's linkage between
two clusters $C$ and $C'$ is the increase in the error sum of squares (or
equivalently, the decrease in variance) when $C$ and $C'$ are merged:
$\delta(C,C') = I(C \cup C') - I(C) - I(C')$, where
$I(C) := \frac{1}{|C|} \sum_{i \in C} \|x_i - \bar{C}\|^2_{\mathbb{R}^d}$ and
$\bar{C} = \frac{1}{n} \sum_{i \in C} x_i$.  It is one of the most widely used
linkages because of its natural interpretation in terms of within/between
cluster variance and because HAC with Ward's linkage can be seen as a greedy
algorithm for least square minimization, similarly to the $k$-means algorithm.
In this paper, the $p$ objects to be clustered are assumed to be ordered by
their indices $i \in \set{1, \dots p}$. We focus on a modification of HAC where
only adjacent clusters are allowed to be merged. This
\emph{adjacency-constrained} HAC is described in Algorithm~\ref{algo:HAC}.

\algnewcommand{\IFor}[1]{\State\algorithmicfor\ #1\ \algorithmicdo}
\algnewcommand{\EndIFor}{\unskip\ \algorithmicend\ \algorithmicfor}
\begin{algorithm}[!h]
  \caption{Adjacency-constrained HAC}
  \begin{algorithmic}[1]
    \State $\cC^0 = (C_i^0)_{1 \leq i \leq p}$ with $C_i^0 = \set{x_i}$
    \Comment{Initalization} \label{algo:HAC1-init}
    \For {$t = 1$ to $p-1$}
      \State $u_t = \arg\min_{u \in \{1, \dots, p-t\}}
      \delta(C^{t-1}_u,C^{t-1}_{u+1})$ \label{algo:HAC-delta}
      \Comment{Best candidate}
      \For {$u = 1$ to $p-t-1$} \label{algo:HAC-update}
        \Comment{Update of $\cC^{t-1}$ into $\cC^{t}$}
        \If{$u < u_t$}
          $C^{t}_{u} = C^{t-1}_{u}$
        \ElsIf{$u = u_t$}
          $C^{t}_{u} = C^{t-1}_u \cup C^{t-1}_{u+1}$
        \ElsIf{$u > u_t$}
          $C^{t}_{u} = C^{t-1}_{u+1}$
        \EndIf
      \EndFor
    \EndFor
  \end{algorithmic}
  \label{algo:HAC}
\end{algorithm}

An implementation in Fortran of this algorithm was provided by
\cite{grimm_CG1987}. This implementation has been integrated in the \RR{}
package \texttt{rioja} \citep{rioja2015}. Its complexity is $\mathcal{O}(p^2)$
(quadratic) both in time and space, which prevents the use of this
implementation for large genomic data sets.

\paragraph{Extension to general similarities.} HAC and adjacency-constrained HAC
are also frequently used when the objects to be clustered do not belong to an
Euclidean space but are described by pairwise dissimilarities. This case has
been formally studied in
\cite{szekely_rizzo_JC2005,strauss_vonmaltitz_PO2017,chavent_etal_CS2018} and
generally involves extending the linkage formula by making an analogy
between the dissimilarity and the Euclidean distance (or the squared Euclidean
distance in some cases). These authors have shown that the simplified update of
the linkage at each step of the algorithm, known as the Lance-Williams formula,
is still valid in this case and that the objective criterion can be interpreted
as the minimization of a so-called ``pseudo inertia''. A similar approach can be
used to extend HAC to data described by an arbitrary similarity between objects,
$S = (s_{ij})_{i,j=1,\ldots,p}$, using a kernel framework as in
\citep{qin_etal_B2003,ahpine_wang_IDA2016}. More precisely, when $S$ is positive
definite, the theory of Reproducing Kernel Hilbert Spaces
\citep{aronszajn_TAMS1950} implies that the data can be embedded in an implicit
Hilbert space. This allows to define Ward's linkage between any two clusters in
terms of the similarity using the so-called ``kernel trick'':
$\forall\, C,\ C' \subset \{1,\ldots,p\}$,
\begin{equation}
  \label{eq:ward-similarity}
  \delta(C,C') = \frac{S(C)}{\abs{C}}  + \frac{S(C')}{\abs{C'}}  -  \frac{ S(C \cup C')}{\abs{C \cup C'}}\,,
\end{equation}
where $S(C) = \sum_{(i,j) \in C^2} s_{ij}$
only depends on $S$ and not on the embedding. Equation~\eqref{eq:ward-similarity} is proved in Section~S1.1 of Supplementary
material.

Extending this approach to the case of a general (that is, possibly
non-positive definite) similarity matrix has been studied in
\cite{miyamoto_etal_SOCPAR2015}. Noting that $(i)$ for a large enough $\lambda$,
the matrix $S_\lambda = S+\lambda I_p$ is positive definite and that $(ii)$
$\delta_{S_\lambda}(C,C') = \delta(C,C') + \lambda$,
\citet[Theorem~1]{miyamoto_etal_SOCPAR2015} concluded that applying Ward's HAC to
$S$ and $S_\lambda$ yields the exact same hierarchy, only shifting the linkage
values by $+\lambda$. This result, which a fortiori holds for the
adjacency-constrained Ward's HAC, justifies the use of
Equation~\eqref{eq:ward-similarity} in the case of a general similarity matrix.

\subsection{Band similarity}
\label{sec:band-similarity}

The implementation provided in \pkg{rioja} takes as an input a $p \times p$
(dense) dissimilarity matrix, making its space complexity quadratic.
Algorithm~\ref{algo:HAC} can be made sub-quadratic in space in situations where
the similarity matrix is sparse (see \cite{ahpine_wang_IDA2016} for similar
considerations in the unconstrained case) or when the similarities can be
computed on the fly, that is, at the time they are required by the algorithm, as
in \cite{dehman_etal_BMCB2015}. However, its time complexity is intrinsically
quadratic in $p$ because all of the $p^2$ similarities are used to compute the
linkages (Algorithm~\ref{algo:HAC}, line \ref{algo:HAC-delta}).

In applications where adjacency-constrained clustering is relevant, such as
Hi-C and GWAS data analysis, this quadratic complexity is a major practical
drawback because $p$ is typically of the order of $10^4$ to $10^5$ for each
chromosome.  Fortunately, in such applications it also makes sense to assume
that the similarity between physically distant objects is small. Specifically,
we assume that $S$ is a band matrix of bandwidth $h+1$, where $h \in \{1\dots
p\}$: $s_{ij} = 0$ for $|i-j| \geq h$.  This assumption is not restrictive,
as it is always fulfilled for $h=p$. However, we will be mostly interested in
the case where $h \ll p$.

In the remainder of this section, we introduce an algorithm with improved time
and space complexity under this band similarity assumption. This algorithm
relies on $(i)$ constant-time calculation of each of the Ward's linkages
involved at line \ref{algo:HAC-delta} of Algorithm~\ref{algo:HAC} using Equation~(\ref{eq:ward-similarity}), and $(ii)$
storage of the candidate fusions in a min-heap. These elements are described in
the next two subsections.

\subsubsection{Ward's linkage as a function of pre-calculated sums}

\begin{figure}[!h]
  \centering
  \includegraphics[width=0.45\linewidth]{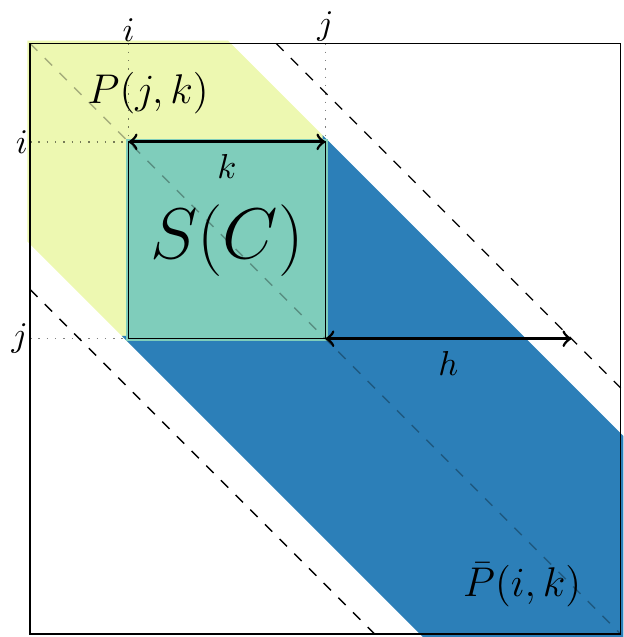}
  \includegraphics[width=0.45\linewidth]{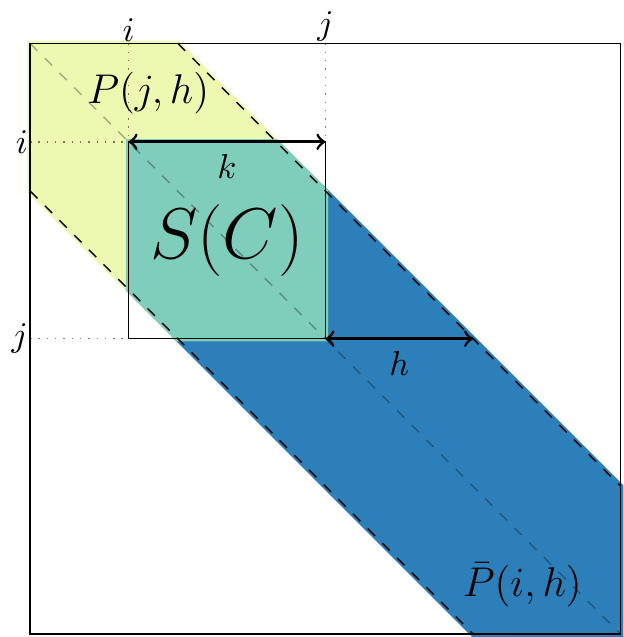}
  \caption{Example of forward pencils (in yellow and green) and backward
pencils (in green and blue), and illustration of Equation~\eqref{eq:pencil-sums}
for cluster $C=\{i, \dots , j-1\}$. Left: cluster smaller than bandwidth ($k
\leq h$); right: cluster larger than bandwidth $k \geq h$.}
  \label{fig:pencil-trick}
\end{figure}

The key point of this subsection is to show that the sums of similarities involved in Equation~(\ref{eq:ward-similarity}) may be expressed as a function of certain pre-calculated sums.
We start by noting that the sum of all similarities in any cluster  $C=\{i,
\dots , j-1\}$ of size $k=j-i$ can easily be obtained from sums of elements in
the first $\min(h,k)$ subdiagonals of $S$. To demonstrate that this is the case we define, for $1 \leq
r,l \leq p$, $P(r,l)$ as the sum of all elements of $S$ in the first $l$
subdiagonals of the upper-right $r \times r$ block of $S$. Formally,
\begin{equation}
  \label{eq::def-pencils}
	P(r,l) = \sum_{1\leq i,j \leq r, \abs{i-j} < l} s_{ij}
\end{equation}
and symmetrically, $\bar{P}(r, l) = P(p+1-r, l)$. This notation is illustrated
in Figure~\ref{fig:pencil-trick}, with $r \in \{i,j\}$. In the left panel, $l=k
\leq h$, while in the right panel, $l=h \leq k$. In both panels, $P(j,
\min(h,k))$ is the sum of elements in the yellow and green regions, while
$\bar{P}(i, \min(h,k))$ is the sum of elements in the green and blue regions.
Because $P$ and $\bar{P}$ are sums of elements in pencil-shaped areas, we call
$P(r,l)$ a \emph{forward pencil} and $\bar{P}(r,l)$ a \emph{backward pencil}.

Figure~\ref{fig:pencil-trick} illustrates that the sum $S_{CC}$ of all
similarities in cluster $C$ can be computed from forward and backward pencils
using the identity:
\begin{equation}
	\label{eq:pencil-sums}
	P(j, h_k) + \bar{P}(i, h_k) = S(C) + P(p, h_k)\,,
\end{equation}
where $h_k:=\min(h,k)$ and $P(p, h_k)$ is the ``full'' pencil of bandwidth
$h_k$
(which also corresponds to $\bar{P}(1, h_k)$). 
The above formula makes it possible to compute $\delta(C,C')$ in constant time
from the pencil sums using Equation~\eqref{eq:ward-similarity}. By construction,
all the bandwidths of the pencils involved are less than $h$. Therefore, only
pencils $P(r,l)$ and $\bar{P}(r,l)$ with $1 \leq r \leq p$ and $1 \leq l \leq h$
have to be pre-computed, so that the total number of pencils to compute and
store is less than $2ph$. These computations can be performed recursively in a
$\mathcal{O}(ph)$ time complexity. Further details about the time and space
complexity of this pencil trick are given in Section 1.2 of the Supplementary
Material.

\subsubsection{Storing candidate fusions in a min-heap}
\label{sec:heapz}

\begin{figure}[!h]
  \centering
  \includegraphics[width=\linewidth]{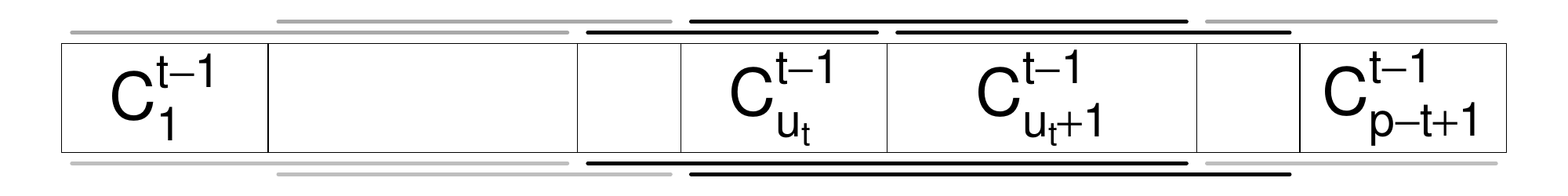}
  \caption{The $t^\textnormal{th}$ merging step in adjacency-constrained HAC in
    Algorithm~\ref{algo:HAC}.  The clusters are represented by rectangular
    cells. Candidate fusions are represented by horizontal bars: above the
    corresponding pair of clusters at step $t$ and below it at step $t+1$,
    assuming that the best fusion is the one between the clusters of indices $u_t$
    and $u_t+1$. Gray bars indicate candidate fusions that are present at both
    steps.}
  \label{fig:algo-chac-merge}
\end{figure}

Iteration $t$ of Algorithm~\ref{algo:HAC} consists in finding the minimum of
$p-t$ elements, corresponding to the candidate fusions between the $p-t+1$
clusters in $\cC^{t-1}$, and merging the corresponding clusters. Storing the
candidate fusions in an \emph{unordered array} and calculating the minimum at
each step would mean a quadratic time complexity. One intuitive strategy would be to make use of the fact that all but 2 to 3 candidate fusions at step $t$ (horizontal bars above
the clusters) are still candidate fusions at step $t-1$, as illustrated by
Figure~\ref{fig:algo-chac-merge}.
However, maintaining a \emph{totally-ordered} list of candidate fusions is not efficient
because the cost of deleting and inserting an element in an ordered list is linear
in $p$, again leading to a quadratic time complexity. Instead, we propose
storing the candidate fusions in a \emph{partially-ordered} data structure called
a \emph{min heap} ~\citep{williams_CACM1964}. This type of structure achieves an
appropriate tradeoff between the cost of maintaining the structure and the cost
of finding the minimum element at each iteration, as illustrated in
Table~\ref{complexityTab} below.
\begin{table}[!h]
  \centering
  \begin{tabular}{r|lll|l}
    & Find $\min$  & Insert & Delete $\min$ & Total\\
    \hline
    Unordered array & $p$ & 1 & $p$ & $p$\\
    Min heap & 1 & $\log(p)$ & $\log(p)$ & $\log(p)$\\
    Ordered array & 1 & $p$ & $p$ & $p$\\
  \end{tabular}
  \caption{Time complexities ($\times \mathcal{O}(1)$) of the three main
elementary operations required by one step of adjacency-constrained clustering
(in columns), for three implementation options (in rows), for a problem of size
$p$.}
  \label{complexityTab}
\end{table}

A min heap is a binary tree such that the value of each node is smaller than
the value of its two children. The advantage of this structure is that all the
operations required in Algorithm~\ref{algo:HAC} to create and maintain the list
of candidate fusions can be done very efficiently. Specifically, at the
beginning of the clustering, the heap is initialized with $p-1$ candidate
fusions in $\mathcal{O}(p \log(p))$. Then, at each of the $p$ iterations is
essentially a combination of a fixed number of elementary operations in at most
$\mathcal{O}(\log(p))$:
\begin{enumerate}
	\item find the best candidate fusion (root of the min heap) in
$\mathcal{O}(1)$;
	\item delete the root of the min heap in $\mathcal{O}(\log(p))$;
	\item insert two possible fusions in the min heap in $\mathcal{O}(\log(p))$
\end{enumerate}

The resulting algorithm is implemented in \pkg{adjclust}. It is described in
details in Algorithm~S2 in Section~S2 of Supplementary Material, where examples
of min heaps are also given.

\subsubsection{Complexity of the proposed algorithm}
\label{sec:complexity-algorithm}

By pre-calculating the $ph$ initial pencils recursively using cumulative sums,
the time complexity of the pre-computation step is $ph$ and the time complexity
of the computation of the linkage of the merged cluster with its two neighbors
is $\mathcal{O}(1)$ (see Section~S1.2 of Supplementary material for further
details). Its total time complexity is thus $\mathcal{O}(p(h + \log(p))$, where
$\mathcal{O}(ph)$ comes from the pre-computation of pencils, and
$\mathcal{O}(p\log(p))$ comes from the $p$ iterations of the algorithm (to merge
clusters from $p$ clusters up to 1 cluster), each of which has a complexity of
$\mathcal{O}(\log(p))$. The space complexity of this algorithm is
$\mathcal{O}(ph)$ because the size of the heap is $\mathcal{O}(p)$ and the space
complexity of the pencil pre-computations is $\mathcal{O}(ph)$. Therefore, the
method achieves a quasi-linear (linearithmic) time complexity and linear space
complexity when $h \ll p$. It is our opinion that this is efficient enough for
analyzing large genomic datasets.

\subsection{Implementation}
\label{sec:choice-number-clusters}

Our method is available in the \RR{} package \pkg{adjclust}, using an
underlying implementation in C and available on
CRAN\footnote{\url{https://cran.r-project.org/package=adjclust}}. Additional features have been implemented to make the package easier to use and results easier to interpret. These include:
\begin{itemize}
  \item plots to display the similarity or dissimilarity together with the
dendrogram and a clustering corresponding to a given level of the hierarchy as
illustrated in Supplementary Figure~S6;
  \item wrappers to use the method with SNP data or Hi-C data that take data
from standard bed files or outputs of the packages \pkg{snpStats} and
\pkg{HiTC} respectively;
  \item a function to guide the user towards a relevant cut of the dendrogram
(and thus a relevant clustering). In practice the
underlying number of clusters is rarely known, and it is important to choose one based on the
data. Two methods are proposed in \pkg{adjclust}: the first is based on a
broken stick model \citep{bennett_NP1996} for the dispersion. Starting from the
root of the dendrogram, the idea is to iteratively check whether the decrease in
within-cluster variance corresponding to the next split can or cannot be explained
by a broken stick model and to stop if it can. To the best of our knowledge
this broken stick strategy is ad hoc in the sense that it does not have a
statistical justification in terms of model selection, estimation of the signal,
or consistency. The second method is based on the slope heuristic that is
statistically justified in the case of segmentation problems
\citep{arlot_etal_p2016,garreau_arlot_p2016}, for which HAC provides an
approximate solution. This later approach is implemented using the
\pkg{capushe} package \citep{arlot_etal_capushe2016}, with a penalty shape
of $p-1 \choose K-1$.
\end{itemize}

Clustering with spatial constraints has many different applications in genomics.
The next two sections illustrate the relevance of our adjacency contraint
clustering approach in dealing with SNP and Hi-C data. In both cases
samples are described by up to a few million variables. All simulations and
figures were performed using the \RR{} package \pkg{adjclust}, version 0.5.7.

\section{Linkage disequilibrium block inference in GWAS}
\label{sec:appli-gwas}

Genome-Wide Association Studies (GWAS) seek to identify causal genomic variants
associated with rare human diseases. The classical statistical approach for
detecting these variants is based on univariate hypothesis testing, with healthy
individuals being tested against affected individuals at each locus.  Given that
an individual's genotype is characterized by millions of SNPs this approach
yields a large multiple testing problem. Due to recombination phenomena, the
hypotheses corresponding to SNPs that are close to each other along the genome
are statistically dependent. A natural way to account for this dependence in the
process is to reduce the number of hypotheses to be tested by grouping and
aggregating SNPs \citep{dehman_etal_BMCB2015,guinot_etal_p2017} based on their
pairwise Linkage Disequilibrium (LD). In particular, a widely used measure of LD
in the context of GWAS is the $r^2$ coefficient, which can be estimated directly
from genotypes measured by genotyping array or sequencing data using standard
methods \citep{snpStats}. The similarity $S = (r^2_{ij})_{i,j}$ induced by LD
can be shown to be a kernel (see Section~S1.3 of Supplementary
material). Identifying blocks of LD may also be useful to define tag SNPs for
subsequent studies, or to characterize the recombination phenomena.

Numerical experiments were performed on a SNP dataset coming from a GWA
study on HIV \citep{dalmasso_etal_PO2008} based on 317k Illumina genotyping
microarrays. For the evaluation we used five data sets corresponding to
five chromosomes that span the typical number of SNPs per chromosome observed on
this array ($p = 23,304$ for chromosome 1, $p=20,811$ for chromosome 6,
$p=14,644$ for chromosome 11, $p=8,965$ for chromosome 16 and $p=5,436$ for
chromosome 21).

 For each dataset, we computed the LD using the function \texttt{ld} of
 \pkg{snpStats}, either for all SNP pairs ($h=p$) or with a reduced number of
 SNP pairs, corresponding to a bandwidth
 $h \in \{100,\ 200,\ 500,\ 1000,\ 2000,\ 5000,\ 10000,\ 20000\}$.  The packages
 \pkg{rioja} \citep{rioja2015} (which requires the full matrix to be given as a
 \texttt{dist} object\footnote{The time needed to compute this matrix was
   50-1000 times larger than the computation of the LD matrix itself. However,
   we did not include this in the total computation time required by \pkg{rioja}
   because we have not tried to optimize it from a computational point of
   view.}) and \pkg{adjclust} with sparse matrices of the class
 \texttt{dgCMatrix} (the default output class of \texttt{ld}) were then used to
 obtain hierarchical clusterings. All simulations were performed on a 64 bit
 Debian 4.9 server, with 512G of RAM, 3GHz CPU (192 processing units) and
 concurrent access.  The available RAM was enough to perform the clustering on
 the full dataset ($h=p$) with \pkg{rioja} although we had previously noticed
 that \pkg{rioja} implementation could not handle more than $8000$ SNPs on a
 standard laptop because of memory issues.

\subsection{Quality of the band approximation}

\begin{figure}[!h]
  \centering
  \includegraphics[width=\linewidth]{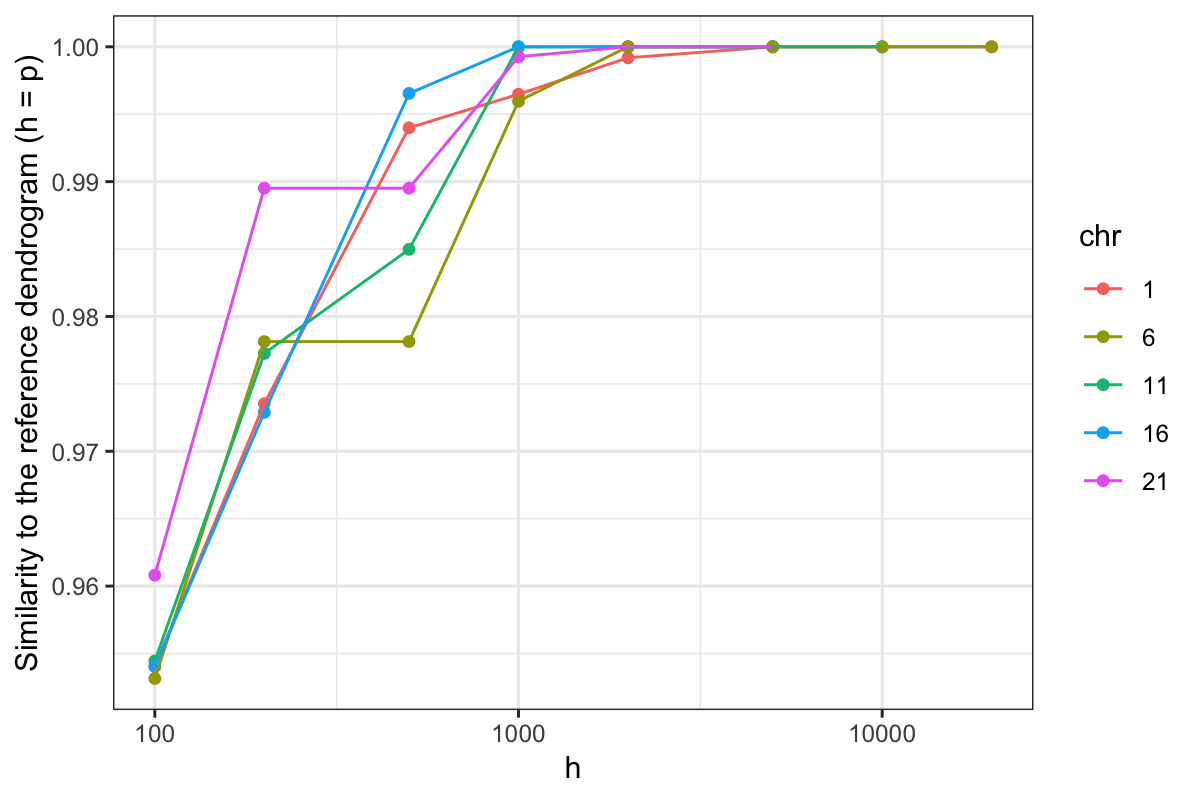}
  \caption{Quality of the band approximation as a function of the bandwidth $h$
    for five different chromosomes.}
  \label{fig:snp-first-diff}
\end{figure}

First, we evaluated the relevance of the band approximation by comparing the
dendrogram obtained with $h < p$ to the reference dendrogram obtained with the
full bandwidth ($h=p$). To perform this comparison we simply recorded the index
$t$ of the last clustering step (among $p-1$) for which all the preceding merges
in the two dendrograms are identical. The quantity $t/(p-1)$ can then be
interpreted as a measure of similarity between dendrograms, ranging from 0 (the
first merges are different) to 1 (the dendrograms are
identical). Figure~\ref{fig:snp-first-diff} displays the evolution of $t/(p-1)$ for
different values of $h$ for the five chromosomes considered here.
For example, for all five chromosomes, at $h = 1000$, the dendrograms differ
from the reference dendrogram only in the last $0.5\%$ of the clustering
step. For $h \geq 2000$ the dendrograms are exactly identical to the reference
dendrogram. We also considered other criteria for evaluating the quality of
the band approximation, including Baker's Gamma correlation coefficient
\citep{Baker1974}, which corresponds to the Spearman correlation between the
ranks of fusion between all pairs of objects. The results obtained with these indices are not shown here because they were consistent with those reported in Figure~\ref{fig:snp-first-diff}.

One important conclusion that may be drawn from these results is that the influence of the bandwidth parameter is the same across chromosomes, that is, across values of $p$  (that range from 5000 to 23000 in this experiment). Therefore, it makes sense to assume that $h$ does not depend on $p$ and that the time and space complexity of our proposed algorithm, which depends on $h$, is indeed quasi-linear in $p$.

\subsection{Scalability and computation times}

\begin{figure}[!h]
  \centering
  \includegraphics[width=\linewidth]{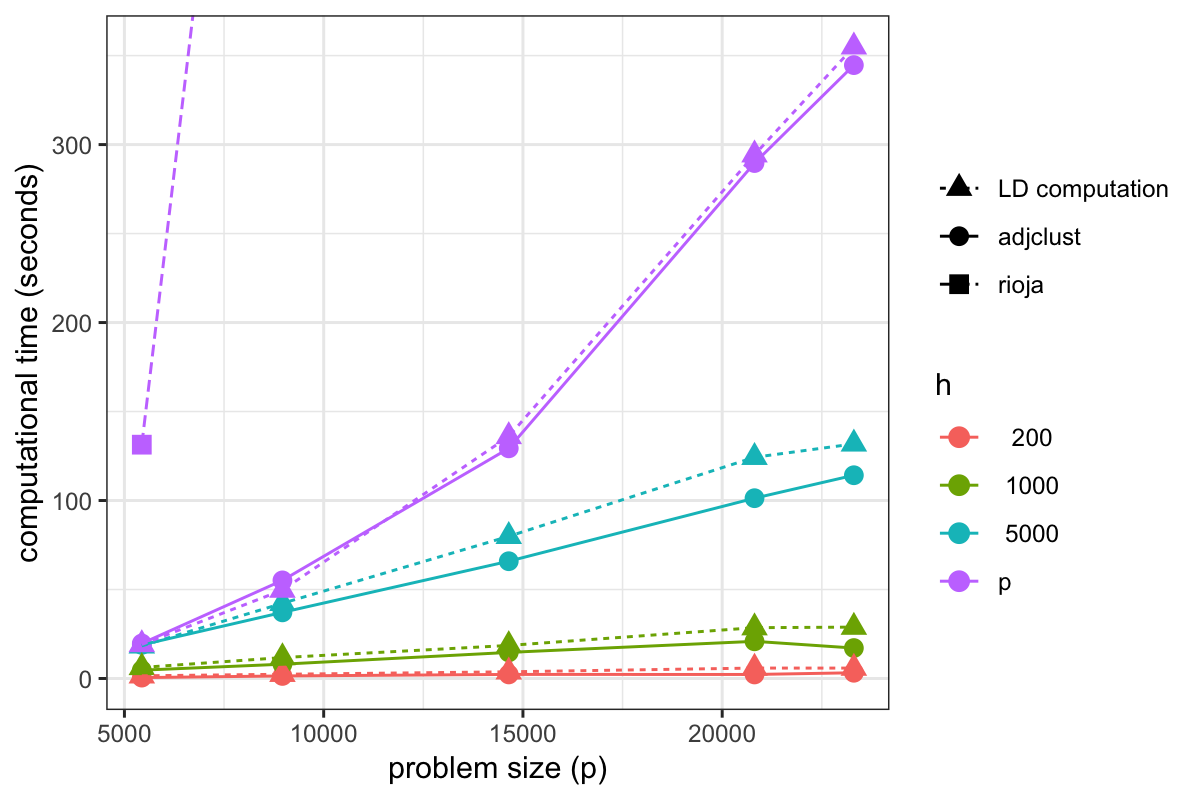}
  \caption{Computation times versus $p$:  LD matrices, for
    CHAC \pkg{rioja} and \pkg{adjclust} with varying
    values for the band $h$.}
  \label{fig::snp-comptime}
\end{figure}

Figure~\ref{fig::snp-comptime} displays the computation time for the LD matrix
(dotted lines) and for the CHAC with respect to the size of the chromosome ($x$
axis), both for \pkg{rioja} (dashed line) and \pkg{adjclust} (solid lines). As
expected, the computation time for \pkg{rioja} did not depend on the bandwidth
$h$, so we only represented $h=p$. For \pkg{adjclust}, the results for varying
bandwidths are represented by different colors. Only the bandwidths
$200$, $1000$, and $5000$ are representend in Figure~\ref{fig::snp-comptime} for
clarity.

Several comments can be made from Figure~\ref{fig::snp-comptime}. First, the
computation times of \pkg{rioja} are much larger than those of \pkg{adjclust},
even when $h=p$ where both methods implement the exact same algorithm. For the
largest chromosome considered here (chromosome 1, $p=23304$), the running time
of \pkg{rioja} is 18900 seconds (more than five hours), compared to 345 seconds
(less than 6 minutes). As expected, the complexity of \pkg{adjclust} with $h=p$
is quadratic in $p$, while it is essentially linear in $p$ for fixed values of $h<p$. For large values of $p$ the gain of the band approximation is substantial: for  $p=23304$ (chromosome 1), the running time of \pkg{adjclust} for $h=1000$ (which is a relevant value in this application according to the results of the preceding section) is of the order of 20 seconds.

We also note that regardless of the value of $h$, the total time needed for the clustering is of the order of (and generally lower than) the time needed for the computation of the LD.

\section{Hi-C analysis}
\label{sec:appli-hi-c}

Hi-C protocol identifies genomic loci that are located nearby in vivo. These
spatial co-locations include intra-chromosomal and inter-chromosomal
interactions. After bioinformatics processing (alignment, filtering, quality
control...), the data are provided as a sparse square matrix with entries that
give the number of reads (contacts) between any given pair of genomic locus bins
at genome scale. Typical sizes of bins are $\sim$40kb, which results in more
than 75,000 bins for the human genome. Constrained clustering or segmentation of
intra-chromosomal maps is a tool frequently used to search for \textit{e.g.},
functional domains (called TADs, Topologically Associating Domains). A number of
methods have been proposed for TAD calling (see \cite{forcato_etal_NM2017} for a
review and comparison), among which the ones proposed by
\cite{fraser_etal_MSB2015,haddad_etal_NAR2017} that take advantage of a
hierarchical clustering, even using a constrained version for the second
reference. In the first article, the authors proceed in two steps with a
segmentation of the data into TADs using a Hidden Markov Model on the
directionality index of Dixon, followed by a greedy clustering on these TADs,
using the mean interaction as a similarity measure between TADs.  Proceeding in
two steps reduces the time required for the clustering, which is $O(p^2)$
otherwise. However, from a statistical and modeling perspective these two steps
would appear redundant. Also, pipelining different procedures (each of them with
their sets of parameters) makes it very difficult to control
errors. \cite{haddad_etal_NAR2017} directly use adjacency-constrained HAC, with
a specific linkage that is not equivalent to Ward's. They do not optimize the
computational time of the whole hierarchy, instead stopping the HAC when a
measure of homogeneity of the cluster created by the last merge falls below a
parameter. Both articles thus highlight the relevance of HAC for exploratory
analysis of Hi-C data. Our proposed approach provides, in addition, a faster way
to obtain an interpretable solution, using the interaction counts as a
similarity and a $h$ similar to the bandwidth of the Dixon index.

\subsection{Data and method}

Data used to illustrate the usefulness of constrained hierarchical clustering
for Hi-C data came from \cite{dixon_etal_N2012,shen_etal_N2012}. Hi-C contact
maps from experiments in mouse embryonic stem cells (mESC), human ESC (hESC),
mouse cortex (mCortex) and human IMR90 Fibroblast (hIMR90) were
downloaded from the authors' website at
\url{http://chromosome.sdsc.edu/mouse/hi-c/download.html} (raw sequence data are
published on the GEO website, accession number GSE35156.

All chromosomes were processed similarly:
\begin{itemize}
	\item counts were $\log$-transformed to reduce the distribution skewness;
	\item constrained hierarchical clustering was computed on $\log$-transformed
data using either the whole matrix ($h=p$) or the sparse approach with a sparse
band size equal to $h = \{0.5p, 0.1p\}$;
	\item model selection was finally performed using both the broken stick
heuristic and the slope heuristic.
\end{itemize}
All computations were performed using the Genotoul cluster.

\subsection{Influence of the bandwidth parameter}

\begin{figure}[ht]
	\centering
  \includegraphics[width=\linewidth]{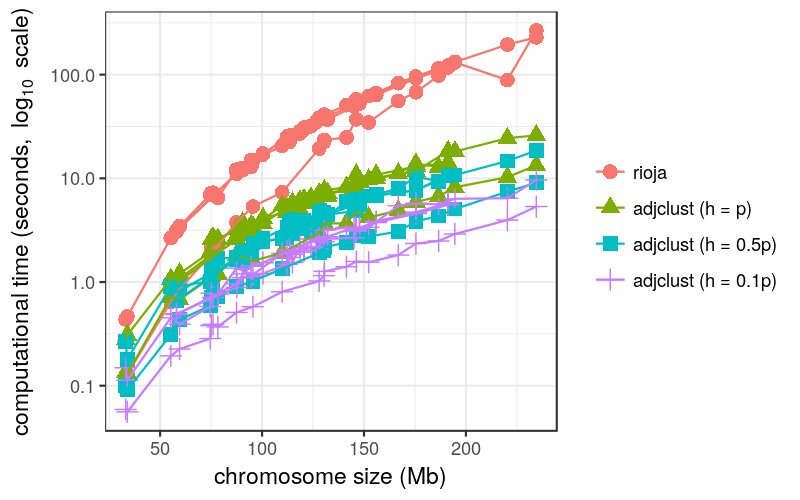}
	\caption{Impact of sparsity on the computational time. Dots that correspond
to the same datasets but different chromosomes are linked by a path.}
	\label{fig::hic-compare-sparse}
\end{figure}

The effect of $h$ (sparse band parameter) on computational time, dendrogram
organization and clustering were assessed. Figure~\ref{fig::hic-compare-sparse}
gives the computational times versus the chromosome size for the three values
of $h$ together with the computational time obtained by the standard version of
constrained hierarchical clustering as implemented in the \RR{} package
\pkg{rioja}.
As expected, the computational time is substantially reduced by the sparse version
(even though not linearly with respect to $h$ because of the preprocessing
step that extracts the band around the diagonal), making the method suitable for
dealing efficiently with a large number of chromosomes and/or a large number of Hi-C
experiments. \pkg{rioja}, that cannot cope efficiently with the sparse band
assumption, requires considerably more computational time (10 times the time needed
by \pkg{adjclust}). In addition, the memory required by the two
approaches is very different: \pkg{adjclust} supports sparse matrix
representation (as implemented in the \RR{} package \pkg{Matrix}), which fits
the way Hi-C matrices are typically stored (usually these matrices are given as
rows with bin number pairs and associated count). For instance, the sparse
version (\texttt{dsCMatrix} class) of the largest chromosome (chromosome 1) in
the hESC data is 23 Mb, as opposed to 231 Mb for the full version. The sparse version of
the smallest chromosome (chromosome 22) is 1.1 Mb, versus 5.2 Mb for the full
version. The sparse version of the $h=0.1p$ band for these two chromosomes is,
respectively, 13.2M and 0.4Mb respectively.

However, this gain in time and space did not impact the results of the method:
the indexes of the first difference were computed between the dendrograms obtained by
the full version ($h=p$) and by the two sparse versions ($h \in \{0.5p,
0.1p\}$) for every chromosome. For most of the clusterings there was no difference in
merge for $h=0.5p$ (with the similarity computed as in
Figure~\ref{fig:snp-first-diff} always larger than 0.9992, and
equal to 1 in more than 3 clusterings out of 4). For $h=0.1p$, the similarity
ranged from 0.9811 to 0.9983. Baker's Gamma index and Rand indices
\citep{hubert_arabie_JC1985} for selected clusterings (both with broken stick
and slope heuristic) confirmed this conclusion (results not shown).

\subsection{Results}

Supplementary Figure~S4 provides the average cluster size for each chromosome
versus the chromosome length. It shows that the average cluster size is fairly
constant among the chromosomes and does not depend on the chromosome length.
Both model selection methods found typical cluster sizes of 1-2 Mb, which is in
line with what is reported in \cite{forcato_etal_NM2017} for some TAD callers.

Supplementary Figure~S5 shows that clusters for a given chromosome (here
chromosome 11 for hIMR90 and chromosome 12 for mCortex) can have different sizes
and also different interpretations: some clusters exhibit a dense interaction
counts (deep yellow) and are thus good TAD candidates whereas a cluster
approximately located between bin 281 and bin 561 in chr12 - mCortex map has
almost no interaction and can be viewed as possibly separating two dense
interaction regions.

\begin{figure}[ht]
	\centering
	\includegraphics[width=\linewidth]{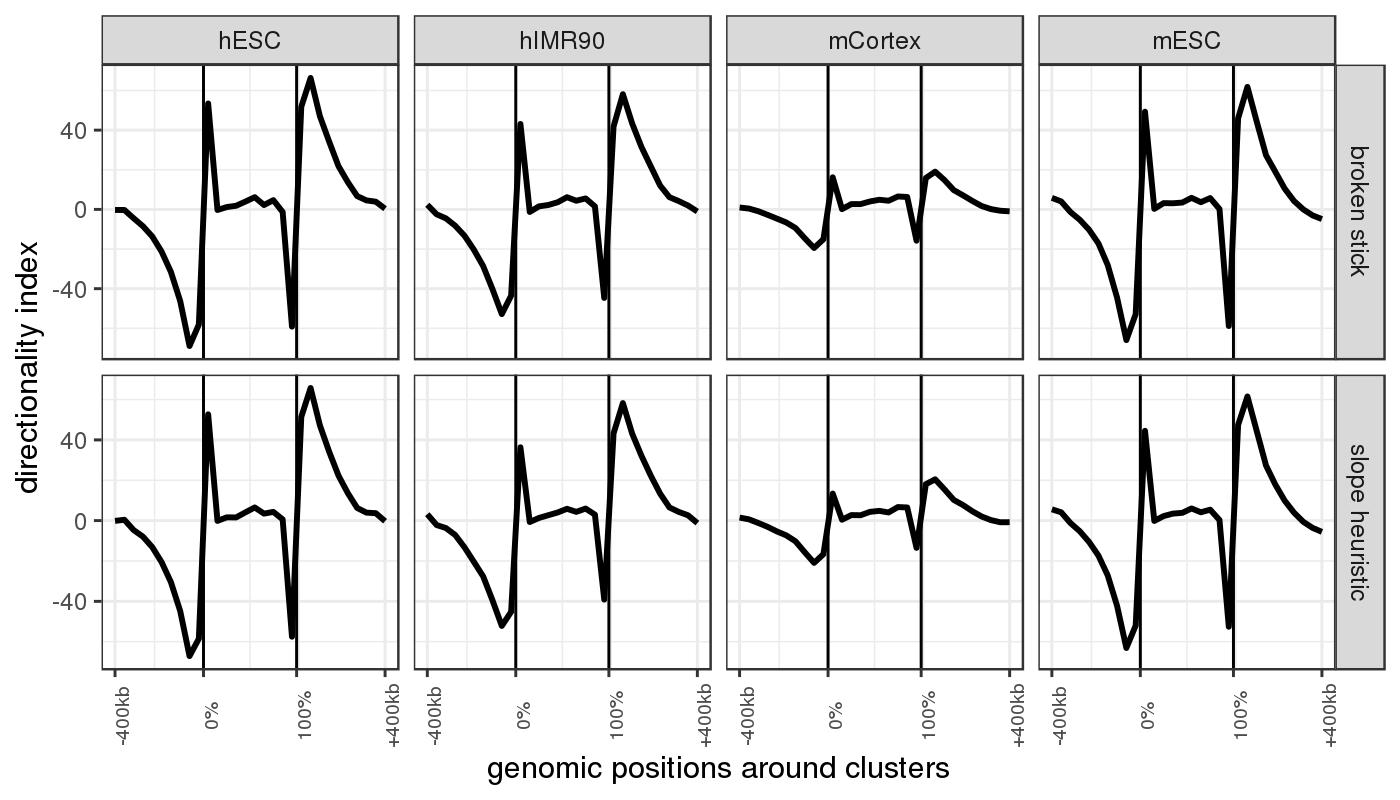}
	\caption{Evolution of the Directionality Index (DI) around clusters (full
version).}
	\label{fig::hic-compare-di}
\end{figure}

The directionality Index (DI, \cite{dixon_etal_N2012}) quantifies a directional
(upstream vs downstream) bias in interaction frequencies, based on a $\chi^2$
statistic. DI is the original method used for TAD calling in Hi-C. Its sign
is expected to change and DI values are expected to show a sharp increase at
TADs boundaries. Figure~\ref{fig::hic-compare-di} displays the average DI, with
respect to the relative bin position within the cluster and the absolute bin
position outside the cluster. The clusters found by constrained HAC show a
relation with DI that is similar to what is expected for standard TADs, with
slightly varying intensities.

Finally, boundaries of TADs are known to be enriched for the insulator binding
protein CTCF \cite{dixon_etal_N2012}. CTCF ChIP-seq peaks were retrieved from
ENCODE \citep{NPC_N2012} and the distribution of the number of the 20\% most
intense peaks was computing at $\pm 400$ Kb of cluster boundaries, as obtained
with the broken stick heuristic (Supplementary Figure~S6). The distribution
also exhibited an enrichment at cluster boundaries, which indicates that the
clustering is relevant with respect to the functional structure of the
chromatin.

\section{Conclusion and prospects}
\label{sec:discussions}

We have proposed an efficient approach to perform constrained hierarchical
clustering based on kernel (or similarity) datasets with several illustrations
of its usefulness for genomic applications. The method is implemented in a
package that is shown to be fast and that currently includes wrappers for
genotyping and Hi-C datasets. The package also provides two possible model
selection procedures to choose a relevant clustering in the hierarchy. The
output of the method is a dendrogram, which can be represented graphically, and
provides a natural hierarchical model for the organization of the objects.

The only tuning parameter in our algorithm is the bandwidth $h$. The numerical
experiments reported in this paper suggest that at least for GWAS and Hi-C
studies, there exists a range of values for $h$ such that $h \ll p$ (which
implies very fast clustering) and the result of the HAC is identical or
extremely close to the clustering obtained for $h=p$. While the range of
relevant values of $h$ will depend on the particular application, an interesting
extension of the present work would be to propose a data-driven choice of $h$ by
running the algorithm on increasing (yet small) values for $h$ on a single
chromosome, and deciding to stop when the dendrogram is stable enough. In
addition, by construction, all groups smaller than $h$ are identical in both
clusterings (with and without the $h$-band approximation).

While HAC is a tool for \emph{exploratory} data analysis, an important
prospect of the present work will be to make use of the low time and memory
footprint of the algorithm in order to perform \emph{inference} on the estimated
hierarchy using stability/resampling-based methods. Such methods could be used
to propose alternative model selection procedures, or to compare hierarchies
corresponding to different biological conditions, which has been shown to be
relevant to Hi-C studies \citep{fraser_etal_MSB2015}.

\section*{Funding}

This work was supported by CNRS project SCALES (Mission ``Osez
l'interdisciplinarité''). The work of G.R. was funded by an ATIGE from
Génopole.

\section*{Acknowledgements}

The authors would like to warmly thank Michel Koskas for very interesting discussions, and for proposing a very elegant alternative implementation.

The authors are grateful to the GenoToul bioinformatics platform (INRA Toulouse,
\url{http://bioinfo.genotoul.fr/}) and its staff for providing computing
facilities.  P. N. and N. V. would like to thank Shubham Chaturvedi for his
contribution to the package \pkg{adjclust} via the R project in google summer of
code 2017.

\bibliography{references}

\begin{thebibliography}{34}
\providecommand{\natexlab}[1]{#1}
\providecommand{\url}[1]{\texttt{#1}}
\expandafter\ifx\csname urlstyle\endcsname\relax
  \providecommand{\doi}[1]{doi: #1}\else
  \providecommand{\doi}{doi: \begingroup \urlstyle{rm}\Url}\fi

\bibitem[Ah-Pine and Wang(2016)]{ahpine_wang_IDA2016}
J.~Ah-Pine and X.~Wang.
\newblock Similarity based hierarchical clustering with an application to text
  collections.
\newblock In H.~Bostr\"om, A.~Knobbe, C.~Soares, and P.~Papapetrou, editors,
  \emph{Proceedings of the 15th International Symposium on Intelligent Data
  Analysis (IDA 2016)}, Lecture Notes in Computer Sciences, pages 320--331,
  Stockholm, Sweden, 2016.
\newblock \doi{10.1007/978-3-319-46349-0}.
\newblock URL \url{https://hal.archives-ouvertes.fr/hal-01437124}.

\bibitem[Arlot et~al.(2016{\natexlab{a}})Arlot, Brault, Baudry, Maugis, and
  Michel]{arlot_etal_capushe2016}
S.~Arlot, V.~Brault, J.-P. Baudry, C.~Maugis, and B.~Michel.
\newblock \emph{capushe: CAlibrating Penalities Using Slope HEuristics},
  2016{\natexlab{a}}.
\newblock URL \url{https://CRAN.R-project.org/package=capushe}.
\newblock R package version 1.1.1.

\bibitem[Arlot et~al.(2016{\natexlab{b}})Arlot, Celisse, and
  Harchaoui]{arlot_etal_p2016}
S.~Arlot, A.~Celisse, and Z.~Harchaoui.
\newblock A kernel multiple change-point algorithm via model selection.
\newblock Preprint arXiv: 1202.3878, 2016{\natexlab{b}}.
\newblock URL \url{https://arxiv.org/abs/1202.3878}.

\bibitem[Aronszajn(1950)]{aronszajn_TAMS1950}
N.~Aronszajn.
\newblock Theory of reproducing kernels.
\newblock \emph{Transactions of the American Mathematical Society}, 68\penalty0
  (3):\penalty0 337--404, 1950.

\bibitem[Baker(1974)]{Baker1974}
F.~B. Baker.
\newblock Stability of two hierarchical grouping techniques case {I}:
  sensitivity to data errors.
\newblock \emph{Journal of the American Statistical Association}, 69\penalty0
  (346):\penalty0 440--445, 1974.
\newblock \doi{10.1080/01621459.1974.10482971}.

\bibitem[Bennett(1996)]{bennett_NP1996}
K.~D. Bennett.
\newblock Determination of the number of zones in a biostratigraphical
  sequence.
\newblock \emph{New Phytologist}, 132\penalty0 (1):\penalty0 155--170, 1996.
\newblock \doi{10.1111/j.1469-8137.1996.tb04521.x}.

\bibitem[Celisse et~al.(2017)Celisse, Marot, Pierre-Jean, and
  Rigaill]{celisse_etal_p2007}
A.~Celisse, G.~Marot, M.~Pierre-Jean, and G.~Rigaill.
\newblock New efficient algorithms for multiple change-point detection with
  kernels.
\newblock Preprint arXiv: 1710.04556, 2017.
\newblock URL \url{https://arxiv.org/abs/1710.04556}.

\bibitem[Chavent et~al.(2018)Chavent, Kuentz-Simonet, Labenne, and
  Saracco]{chavent_etal_CS2018}
M.~Chavent, V.~Kuentz-Simonet, A.~Labenne, and J.~Saracco.
\newblock {ClustGeo2}: an {R} package for hierarchical clustering with spatial
  constraints.
\newblock \emph{Computational Statistics}, 33\penalty0 (4):\penalty0
  1799--1822, 2018.
\newblock \doi{10.1007/s00180-018-0791-1}.

\bibitem[Clayton(2015)]{snpStats}
D.~Clayton.
\newblock \emph{{snpStats}: {SnpMatrix} and {XSnpMatrix} classes and methods},
  2015.
\newblock R package version 1.24.0.

\bibitem[Dalmasso et~al.(2008)Dalmasso, Carpentier, Meyer, Rouzioux, Goujard,
  Chaix, Lambotte, Avettand-Fenoel, Le~Clerc, de~Senneville, Deveau, Boufassa,
  Debr\'e, Delfraissy, Broet, and Theodorou]{dalmasso_etal_PO2008}
C.~Dalmasso, W.~Carpentier, L.~Meyer, C.~Rouzioux, C.~Goujard, M.-L. Chaix,
  O.~Lambotte, V.~Avettand-Fenoel, S.~Le~Clerc, L.~D. de~Senneville, C.~Deveau,
  F.~Boufassa, P.~Debr\'e, J.-F. Delfraissy, P.~Broet, and I.~Theodorou.
\newblock {D}istinct genetic loci control plasma {HIV-RNA} and cellular
  {HIV-DNA} levels in {HIV-1} infection: the {ANRS} {G}enome {W}ide
  {A}ssociation 01 study.
\newblock \emph{PLoS ONE}, 3\penalty0 (12):\penalty0 e3907, 2008.
\newblock \doi{10.1371/journal.pone.0003907}.

\bibitem[Dehman et~al.(2015)Dehman, Ambroise, and
  Neuvial]{dehman_etal_BMCB2015}
A.~Dehman, C.~Ambroise, and P.~Neuvial.
\newblock Performance of a blockwise approach in variable selection using
  linkage disequilibrium information.
\newblock \emph{BMC Bioinformatics}, 16\penalty0 (1):\penalty0 148, 2015.
\newblock \doi{10.1186/s12859-015-0556-6}.

\bibitem[Dixon et~al.(2012)Dixon, Selvaraj, Yue, Kim, Li, Shen, Hu, Liu, and
  Ren]{dixon_etal_N2012}
J.~Dixon, S.~Selvaraj, F.~Yue, A.~Kim, Y.~Li, Y.~Shen, M.~Hu, J.~Liu, and
  B.~Ren.
\newblock Topological domains in mammalian genomes identified by analysis of
  chromatin interactions.
\newblock \emph{Nature}, 485:\penalty0 376--380, 2012.
\newblock \doi{10.1038/nature11082}.

\bibitem[{ENCODE Project Consortium}(2012)]{NPC_N2012}
{ENCODE Project Consortium}.
\newblock An integrated encyclopedia of {DNA} elements in the human genome.
\newblock \emph{Nature}, 489:\penalty0 57--74, 2012.
\newblock \doi{10.1038/nature11247}.

\bibitem[Forcato et~al.(2017)Forcato, Nicoletti, Pal, Livi, Ferrari, and
  Bicciato]{forcato_etal_NM2017}
M.~Forcato, C.~Nicoletti, K.~Pal, C.~Livi, F.~Ferrari, and S.~Bicciato.
\newblock Comparison of computational methods for {Hi-C} data analysis.
\newblock \emph{Nature Methods}, 14\penalty0 (7):\penalty0 679--685, 2017.

\bibitem[Fraser et~al.(2015)Fraser, Ferrai, Chiariello, Schueler, Rito,
  Laudanno, Barbieri, Moore, Kraemer, Aitken, Xie, Morris, Itoh, Kawaji,
  Jaeger, Hayashizaki, Carninci, Forrest, {The FANTOM Consortium}, Semple,
  Dostie, Pombo, and Nicodemi]{fraser_etal_MSB2015}
J.~Fraser, C.~Ferrai, A.~Chiariello, M.~Schueler, T.~Rito, G.~Laudanno,
  M.~Barbieri, B.~Moore, D.~Kraemer, S.~Aitken, S.~Xie, K.~Morris, M.~Itoh,
  H.~Kawaji, I.~Jaeger, Y.~Hayashizaki, P.~Carninci, A.~Forrest, {The FANTOM
  Consortium}, C.~Semple, J.~Dostie, A.~Pombo, and M.~Nicodemi.
\newblock Hierarchical folding and reorganization of chromosomes are linked to
  transcriptional changes in cellular differentiation.
\newblock \emph{Molecular Systems Biology}, 11:\penalty0 852, 2015.
\newblock \doi{10.15252/msb.20156492}.

\bibitem[Gabriel et~al.(2002)Gabriel, Schaffner, Nguyen, Moore, Roy,
  Blumenstiel, Higgins, DeFelice, Lochner, Faggart, Liu-Cordero, Rotimi,
  Adeyemo, Cooper, Ward, Lander, Daly, and Altshuler]{gabriel_etal_S2002}
S.~B. Gabriel, S.~F. Schaffner, H.~Nguyen, J.~M. Moore, J.~Roy, B.~Blumenstiel,
  J.~Higgins, M.~DeFelice, A.~Lochner, M.~Faggart, S.~N. Liu-Cordero,
  C.~Rotimi, A.~Adeyemo, R.~Cooper, R.~Ward, E.~S. Lander, M.~J. Daly, and
  D.~Altshuler.
\newblock The structure of haplotype blocks in the human genome.
\newblock \emph{Science}, 296\penalty0 (5576):\penalty0 2225--2229, 2002.
\newblock \doi{10.1126/science.1069424}.

\bibitem[Garreau and Arlot(2016)]{garreau_arlot_p2016}
D.~Garreau and S.~Arlot.
\newblock Consistent change-point detection with kernels.
\newblock arXiv preprint arXiv:1612.04740, 2016.
\newblock URL \url{https://arxiv.org/abs/1612.04740}.

\bibitem[Grimm(1987)]{grimm_CG1987}
E.~Grimm.
\newblock {CONISS}: a fortran 77 program for stratigraphically constrained
  analysis by the method of incremental sum of squares.
\newblock \emph{Computers \& Geosciences}, 13\penalty0 (1):\penalty0 13--35,
  1987.

\bibitem[Guinot et~al.(2017)Guinot, Szafranski, Ambroise, and
  Samson]{guinot_etal_p2017}
F.~Guinot, M.~Szafranski, C.~Ambroise, and F.~Samson.
\newblock Learning the optimal scale for {GWAS} through hierarchical {SNP}
  aggregation.
\newblock Preprint arXiv: 1710.01085, 2017.
\newblock URL \url{https://arxiv.org/abs/1710.01085}.

\bibitem[Haddad et~al.(2017)Haddad, Vaillant, and Jost]{haddad_etal_NAR2017}
N.~Haddad, C.~Vaillant, and D.~Jost.
\newblock {IC}-{F}inder: inferring robustly the hierarchical organization of
  chromatin folding.
\newblock \emph{Nucleic Acids Research}, 45\penalty0 (10):\penalty0 e81, 2017.
\newblock \doi{10.1093/nar/gkx036}.

\bibitem[Harchaoui and Capp{\'e}(2007)]{harchaoui_cape_WSSP2007}
Z.~Harchaoui and O.~Capp{\'e}.
\newblock Retrospective mutiple change-point estimation with kernels.
\newblock In \emph{Proceedings of the 14th Workshop on Statistical Signal
  Processing (SSP'07)}, pages 768--772, Madison, WI, USA, 2007. IEEE.
\newblock \doi{10.1109/SSP.2007.4301363}.

\bibitem[Hocking et~al.(2013)Hocking, Schleiermacher, Janoueix-Lerosey, Boeva,
  Cappo, Delattre, Bach, and Vert]{hocking_etal_BMCB2013}
T.~D. Hocking, G.~Schleiermacher, I.~Janoueix-Lerosey, V.~Boeva, J.~Cappo,
  O.~Delattre, F.~Bach, and J.-P. Vert.
\newblock Learning smoothing models of copy number profiles using breakpoint
  annotations.
\newblock \emph{BMC Bioinformatics}, 14\penalty0 (1):\penalty0 164, 2013.
\newblock \doi{10.1186/1471-2105-14-164}.

\bibitem[Hubert and Arabie(1985)]{hubert_arabie_JC1985}
L.~Hubert and P.~Arabie.
\newblock Comparing partitions.
\newblock \emph{Journal of Classification}, 2\penalty0 (1):\penalty0 193--218,
  1985.
\newblock \doi{10.1007/BF01908075}.

\bibitem[Juggins(2018)]{rioja2015}
S.~Juggins.
\newblock \emph{rioja: Analysis of Quaternary Science Data}, 2018.
\newblock URL \url{https://cran.r-project.org/package=rioja}.
\newblock R package version 0.9-15.1.

\bibitem[Kaufman and Rousseeuw(2009)]{kaufman_rousseeuw_FGDICA2009}
L.~Kaufman and P.~J. Rousseeuw.
\newblock \emph{Finding Groups in Data: an Introduction to Cluster Analysis},
  volume 344 of \emph{Wiley Series in Probability and Statistics}.
\newblock John Wiley \& Sons, Hoboken, NJ, USA, 2009.
\newblock ISBN 9780471878766.
\newblock \doi{10.1002/9780470316801}.

\bibitem[Lebart(1978)]{lebart_CAD1978}
L.~Lebart.
\newblock Programme d'agr{\'e}gation avec contraintes.
\newblock \emph{Les Cahiers de l'Analyse des Donn{\'e}es}, 3\penalty0
  (3):\penalty0 275--287, 1978.
\newblock URL \url{http://www.numdam.org/item?id=CAD_1978__3_3_275_0}.

\bibitem[Miyamoto et~al.(2015)Miyamoto, Abe, Endo, and
  Takeshita]{miyamoto_etal_SOCPAR2015}
S.~Miyamoto, R.~Abe, Y.~Endo, and J.~Takeshita.
\newblock Ward method of hierarchical clustering for non-{E}uclidean similarity
  measures.
\newblock In \emph{Proceedings of the VIIth International Conference of Soft
  Computing and Pattern Recognition (SoCPaR 2015)}, 2015.

\bibitem[Picard et~al.(2005)Picard, Robin, Lavielle, Vaisse, and
  Daudin]{picard_etal_BMCB2005}
F.~Picard, S.~Robin, M.~Lavielle, C.~Vaisse, and J.-J. Daudin.
\newblock A statistical approach for array-{CGH} data analysis.
\newblock \emph{BMC Bioinformatics}, 6\penalty0 (27):\penalty0 1471--2105,
  2005.
\newblock \doi{10.1186/1471-2105-6-27}.

\bibitem[Qin et~al.(2003)Qin, Lewis, and Noble]{qin_etal_B2003}
J.~Qin, D.~P. Lewis, and W.~S. Noble.
\newblock Kernel hierarchical gene clustering from microarray expression data.
\newblock \emph{Bioinformatics}, 19\penalty0 (16):\penalty0 2097--2104, 2003.
\newblock \doi{10.1093/bioinformatics/btg288}.

\bibitem[Reuter et~al.(2015)Reuter, Spacek, and Snyder]{reuter_etal_MC2015}
J.~A. Reuter, D.~V. Spacek, and M.~P. Snyder.
\newblock High-throughput sequencing technologies.
\newblock \emph{Molecular Cell}, 58\penalty0 (4):\penalty0 586--597, 2015.
\newblock \doi{10.1016/j.molcel.2015.05.004}.

\bibitem[Shen et~al.(2012)Shen, Yu, McCleary, Ye, Edsall, Kuan, Wagner, Dixon,
  Lee, Lobanenkov, and Ren]{shen_etal_N2012}
Y.~Shen, F.~Yu, D.~F. McCleary, Z.~Ye, L.~Edsall, S.~Kuan, U.~Wagner, J.~Dixon,
  L.~Lee, V.~V. Lobanenkov, and B.~Ren.
\newblock A map of the \textit{cis}-regularoty sequence in the mouse genome.
\newblock \emph{Nature}, 488:\penalty0 116--120, 2012.
\newblock \doi{10.1038/nature11243}.

\bibitem[Strauss and von Maltitz(2017)]{strauss_vonmaltitz_PO2017}
T.~Strauss and M.~J. von Maltitz.
\newblock Generalising {W}ard's method for use with {M}anhattan distances.
\newblock \emph{PLoS ONE}, 12:\penalty0 e0168288, 2017.
\newblock \doi{10.1371/journal.pone.0168288}.

\bibitem[Sz\'ekely and Rizzo(2005)]{szekely_rizzo_JC2005}
G.~J. Sz\'ekely and M.~L. Rizzo.
\newblock Hierarchical clustering via joint between-within distances: extending
  {W}ard's minimum variance method.
\newblock \emph{Journal of Classification}, 22\penalty0 (2):\penalty0 151--183,
  2005.
\newblock \doi{10.1007/s00357-005-0012-9}.

\bibitem[Williams(1964)]{williams_CACM1964}
J.~W.~J. Williams.
\newblock Algorithm 232 - heapsort.
\newblock \emph{Communications of the {ACM}}, 7\penalty0 (6):\penalty0
  347--348, 1964.
\newblock \doi{10.1145/512274.512284}.

\end{thebibliography}



\section*{Supplementary material}

# Supplementary file for "Adjacency-constrained hierarchical clustering of a band similarity matrix with application to Genomics"

Christophe Ambroise[1], Alia Dehman[2], Pierre Neuvial[3],
Guillem Rigaill[1,4] and Nathalie Vialaneix[5]

[1] Laboratoire de Mathématiques et Modélisation d'Evry, UMR CNRS 8071, Université d'Evry Val d'Essonne, 23 boulevard de France, 91037 Evry, France.
[2] Hyphen-stat, 195 Route d'Espagne, 31036 Toulouse, France.
[3] Institut de Mathématiques de Toulouse, UMR5219 CNRS, Université de Toulouse, UPS IMT, F-31062 Toulouse Cedex 9, France.
[4] Institute of Plant Sciences Paris Saclay IPS2, CNRS, INRA, Gif sur Yvette, France.
[5] MIAT, Université de Toulouse, INRA, Castanet-Tolosan, France.

## Contents

# S1 Supplementary methods

## S1.1 Proof of Equation (1)

*Proof of Equation (1).* The theory of Reproducing Kernel Hilbert Spaces Aronszajn [1950] makes is possible to generalize the definition of Ward-based HAC to the case where the similarity matrix $S$ describing the objects to cluster is a kernel, *i.e.*, a positive definite symmetric matrix. In this case, there exists a unique Hilbert space $(\mathcal{H}, \langle .,. \rangle_\mathcal{H})$ and a feature map $\phi: \mathcal{X} \to \mathcal{H}$, where $\mathcal{X}$ denotes the arbitrary set in which the objects, $\{x_1, \ldots, x_p\}$, described by $s$ are defined, such that the kernel $s$ corresponds to the dot product in $\mathcal{H}$:

$$s_{ij} = \langle \phi(x_i), \phi(x_j) \rangle_\mathcal{H}. \tag{S1}$$

Following Murtagh and Legendre [2014], since the feature space $\mathcal{H}$ is Euclidean, Ward's linkage may be written as

$$\forall C, C' \subset \{x_1, \ldots, x_p\},\ C \cap C' = \varnothing, \qquad \delta(C, C') = \frac{|C||C'|}{|C| + |C'|} \left\| \bar{C} - \bar{C}' \right\|_\mathcal{H}^2,$$

1

where for any cluster $C$, $\bar{C} := \frac{1}{|C|}\sum_{i \in C} \phi(x_i)$ is the center of gravity of $C$ in $\mathcal{H}$ and $\|\cdot\|_{\mathcal{H}}^2$ is the norm associated to the scalar product in $\mathcal{H}$. However,

$$\begin{aligned}\|\bar{C} - \bar{C}'\|_{\mathcal{H}}^2 &= \langle \bar{C} - \bar{C}', \bar{C} - \bar{C}'\rangle_{\mathcal{H}} \\ &= \langle \bar{C}, \bar{C}\rangle_{\mathcal{H}} + \langle \bar{C}', \bar{C}'\rangle_{\mathcal{H}} - 2\langle \bar{C}, \bar{C}'\rangle_{\mathcal{H}}.\end{aligned}$$

Then, Equation (S1) yields $\forall C, C' \subset \{x_1, \ldots, x_p\}$,

$$\langle \bar{C}, \bar{C}'\rangle_{\mathcal{H}} = \frac{1}{|C||C'|}\langle \sum_{i \in C}\phi(x_i), \sum_{i \in C'}\phi(x_i)\rangle_{\mathcal{H}} = \frac{1}{|C||C'|}S_{CC'}.$$

where $S_{CC'} = \sum_{i \in C, j \in C'} s_{ij}$. This implies

$$\delta(C, C') = \frac{|C||C'|}{|C| + |C'|}\left(\frac{S_{CC}}{|C|^2} + \frac{S_{C'C'}}{|C'|^2} - \frac{2S_{CC'}}{|C||C'|}\right)$$

where $S_{CC} = S(C)$ and $S_{C'C'} = S(C')$. To conclude, we note that

$$\begin{aligned}S(C \cup C') &= S(C) + S(C') + 2\sum_{i \in C, j \in C'} s_{ij} \\ &= S_{CC} + S_{C'C'} + 2S_{CC'},\end{aligned}$$

so that

$$\begin{aligned}\delta(C, C') &= \frac{|C||C'|}{|C| + |C'|}\left(\left(\frac{1}{|C|^2} + \frac{1}{|C||C'|}\right)S(C) + \left(\frac{1}{|C'|^2} + \frac{1}{|C||C'|}\right)S(C') - \frac{1}{|C||C'|}S(C \cup C')\right) \\ &= \frac{S(C)}{|C|} + \frac{S(C')}{|C'|} - \frac{S(C \cup C')}{|C \cup C'|}\end{aligned}$$

which concludes the proof. $\square$

Because $\delta(C, C')$ is explicitely written in terms of $S$ only, Ward's HAC can be performed implicitly in the feature space, without an explicit calculation. In particular, the mapping $\phi$ itself does not need to be known explicitly. This property is known as the kernel trick. In our paper, the kernel trick is used to write the distance between the centers of gravity in the feature space as a function of $S$ only. Ward's HAC is then obtained by iteratively updating the matrix of similarities between all pairs of centers of gravity after each successive merge. To the best of our knowledge, the formulation of Ward's linkage in terms of sums of elements of the similarity matrix $S$ has never been explicitly written in the form of Equation (1) even if kernel-based HAC has already been proposed by Qin et al. [2003], Ah-Pine and Wang [2016].

### S1.2 Time and space complexity of the pencil trick

**Space complexity** The pencil trick is based on the pre-computation of backward and forward pencils as defined in Equation (2). By construction, the number of all the bandwidths of the pencils involved is less than $h$. Therefore, only pencils $P(r,l)$ and $\bar{P}(r,l)$ with $1 \leq r \leq p$ and $1 \leq l \leq h$ have to be pre-computed and the total number of pencils to compute and stored is less than $2ph$. The space complexity is thus $\mathcal{O}(ph)$ for pencils.

**Time complexity** The contribution of the pencil trick to the algorithm complexity is divided into:

- during the initialization of the algorithm, the **pre-computation of the backward and forward pencils**. This can be done efficiently by a recursive computation, as described in Algorithm S1. This algorithm is based on the computation of less than $ph$ quantities, all having a complexity equal to 1. The total time complexity of the method is thus $\mathcal{O}(ph)$;

---
**Algorithm S1** Pencil trick: Precomputation of pencils by a recursive strategy
---
1: **for** $i = 1$ to $p$ **do**
2:     $P(i, 1) \leftarrow s(i, i)$                                                                                            ▷ $\mathcal{O}(1)$ for every $i$
3:     **for** $l = 2$ to $\min(h, p + 1 - i)$ **do**
4:         $P(i, l) \leftarrow P(i, l - 1) + s(i, i + l - 1)$             ▷ $\mathcal{O}(1)$ for every $i$ and every $l$
5:     **end for**
6: **end for**

---

- during the call to HEAP.INSERT (see Section S2 and Algorithm S2), the computation of the linkage values between the new(ly merged) cluster and its right and left neighbors. The time complexity of each linkage computation is $\mathcal{O}(1)$ (constant) because according to Equations (1) and (4), $\delta$ is a function of a constant number of pencils.

## S1.3 Linkage disequilibrium and kernel

If SNP values at position $i$ are modeled by a binary random variable $Z_i$ which is the indicator of the presence of minor allele for this locus, it is standard to make the asumption that

$$Z_i \sim \mathcal{B}(p_i).$$

The linkage disequilibrium (LD) between locus $i$ and locus $j$ is defined as the covariance between the two corresponding random variables:

$$D_{ij} = p_{ij} - p_i p_j = E[Z_i Z_j] - p_i p_j = \text{Cov}(Z_i, Z_j).$$

For practical use the measure is normalized to be between zero and one. Two popular choices of normalization are the squared correlation

$$r^2(i, j) = \text{Cor}(Z_i, Z_j)^2$$

and

$$d'_{ij}(i, j) = \frac{D(i, j)}{\max_{ij} D(i, j)}.$$

In this paper, we consider the $r^2$ which is a classical choice in the context of association studies. More precisely, given observations of the genotypes, or $n$ individuals, we denote by $\boldsymbol{z}_i$ the $2n$-dimensional vector of normalized allele values of locus $i$ for the $2n$ genotypes of the $n$ individuals and estimates the LD with $k(i, j) := \left( \sum_{\ell=1}^{2n} \boldsymbol{z}_{i\ell} \boldsymbol{z}_{j\ell} \right)^2$.

It is possible to prove that this estimation, $k$, is a kernel, *i.e.*, a positive definite symmetric matrix. This result comes from the fact that $k(i, j)$ can be re-written as:

$$\sum_{\ell, \ell'=1}^{2n} (z_{i\ell} z_{j\ell})(z_{i\ell'} z_{j\ell'}).$$

Defining the mapping $\Phi$ from $\mathbb{R}^{(2n)}$ to $\mathbb{R}^{(2n)^2}$ such that

$$\forall (\ell, \ell') \in \{1, \ldots, 2n\}^2, \qquad \Phi(\boldsymbol{z})_{2n(\ell-1)+\ell'} = z_\ell z_{\ell'},$$

we have:

$$k(i,j) = \langle \boldsymbol{z}_i, \boldsymbol{z}_j \rangle^2 = \langle \Phi(\boldsymbol{z}_i), \Phi(\boldsymbol{z}_j) \rangle,$$

which concludes the proof since $k$ is expressed as a dot product with the feature map $i \to \Phi(\boldsymbol{z}_i)$.

Notice that, when the available data are SNPs, computing the square correlation between of two loci raises the additional problem of unknown haplotype phase. Indeed, with association study, we observe locus values for pairs of chromosomes and not for specific chromosomes. For each locus, SNP data give access to a $3 \times 3$ contingency table of genotype counts for each pair of loci, while we would like to have the $2 \times 2$ table of diplotype counts. Nevertheless, these are classical approaches for estimating the diplotypes from the genotypes.

## S2 Algorithm

In this section we provide a detailed description of the method presented in Section 2.2.2 of the article, which is implemented in the **adjclust** package. We also give illustrations of the first steps of this algorithm when applied to the RLGH data set provided in the package **rioja**, which data are relative abundances of 41 taxa in $p = 20$ stratigraphic samples. A detailed description of this data set is provided in the help of the RLGH data set.

While the method in Section 2.2.2 is described in terms of clusters, Algorithm S2 is best expressed in terms of candidate fusions. The initialization step (lines 1 to 3) consists in building the heap of $p - 1$ candidate fusions between the $p$ adjacent items. At the end of this step, the root of the heap contains the best fusion among such fusions. This is illustrated in Figure S1 for the RLGH data set. The best candidate fusion, which is by definition the root of the tree, consists in merging $\{4\}$ and $\{5\}$. It is highlighted in violet and the two "neighbor fusions", *i.e.*, the fusions that involve either $\{4\}$ or $\{5\}$, are highlighted in pink. The initialization step has a $\mathcal{O}(p \log(p))$ time complexity because the complexity of inserting each of the $p-1$ elements in the heap is upper bounded by the maximal depth of the heap, that is, $\log_2(p)$.

As stated in Section 2.2.3, the merging step consists in finding the best candidate fusion (line 5), removing it from the heap (line 6) and inserting (up to) two possible fusions (lines 11-12). The other lines of the algorithm explain how the information regarding the adjacent fusions and clusters are retrieved and updated. The notation is illustrated in Figure S2, elaborating on the example of Figure 1 of the article.

The state of the heap after the first fusion is illustrated by Figure S3, where the two new candidate fusions are highlighted in yellow. The two fusions highlighted in grey are the neighbors of the first fusion.

In Algorithm S2 we have omitted several points for simplicity and conciseness of exposition. For a more complete description, the following remarks can be made:

1. The calculation of the linkage is not mentioned explicitly in the calls to `Heap.Insert`. As explained in the main text and detailed in Section S1.2, the linkage between any two clusters can be calculated (in constant time) from pre-calculated pencil sums.

2. Algorithm S2 should take appropriate care of cases when the best fusion involves the first or last cluster. In particular, only one new fusion is defined and inserted in such cases. This is taken care of in the **adjclust** package, but not in Algorithm S2 for simplicity of exposition.

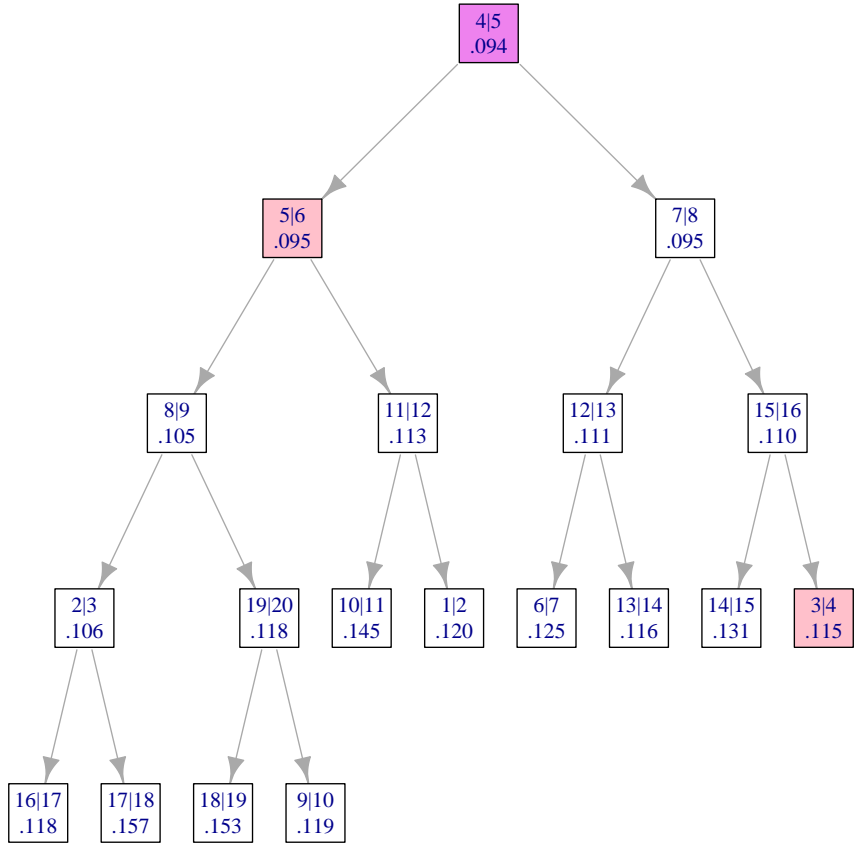

**Figure S1:** Min heap after the initialization step of the RLGH data set. Each node corresponds to a candidate fusion, and is represented by a label of the form $i|i+1$ giving the indices of the items to be merged, and (ii) the value of the corresponding linkage $\delta(\{i\}, \{j\})$. The nodes corresponding to the best fusion and the two neighbor fusions are highlighted.

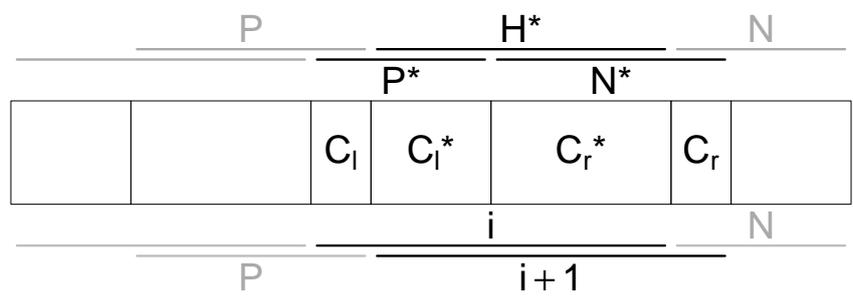

**Figure S2:** Illustration of the result of a merging step in Algorithm S2.

5

**Algorithm S2** adjclust: Adjacency-constrained Ward's HAC of a band similarity

1: **for** $i = 1$ to $p - 1$ **do**  ▷ Initialization
2:     HEAP.INSERT(id $= i$, left $= \{i\}$, right $= \{i+1\}$, prev$=i-1$, next$=i+1$)
3: **end for**
4: **for** $i = p$ to $2p - 1$ by 2 **do**  ▷ Merging
5:     $H^\star \leftarrow$ HEAP.GETROOT( )  ▷ Find best fusion
6:     HEAP.DELETEROOT( )  ▷ Delete min element
7:     $P^\star \leftarrow$ PREV($H^\star$); $N^\star \leftarrow$ NEXT($H^\star$)  ▷ Preceding/next fusion
8:     $P \leftarrow$ PREV($P^\star$); $N \leftarrow$ NEXT($N^\star$)  ▷ Preceding/next fusion
9:     $C_l \leftarrow$ LEFT($P^\star$); $C_{l^\star} \leftarrow$ RIGHT($P^\star$)  ▷ Corresponding clusters
10:    $C_{r^\star} \leftarrow$ LEFT($N^\star$); $C_r \leftarrow$ RIGHT($N^\star$)  ▷ Corresponding clusters
11:    HEAP.INSERT(id$=i$, left$=C_l$, right$=C_{l^\star} \cup C_{r^\star}$, prev$=P$, next$=H^\star$)  ▷ Add new fusion
12:    HEAP.INSERT(id$=i+1$, left$=C_{l^\star} \cup C_{r^\star}$, right$=C_r$, prev$=H^\star$, next$=N$)  ▷ Add new fusion
13:    TAG($P^\star$); TAG($N^\star$)  ▷ Tag inactive fusions
14:    NEXT($P$) $\leftarrow i$; NEXT($N$) $\leftarrow i+1$  ▷ Update neighbors
15: **end for**

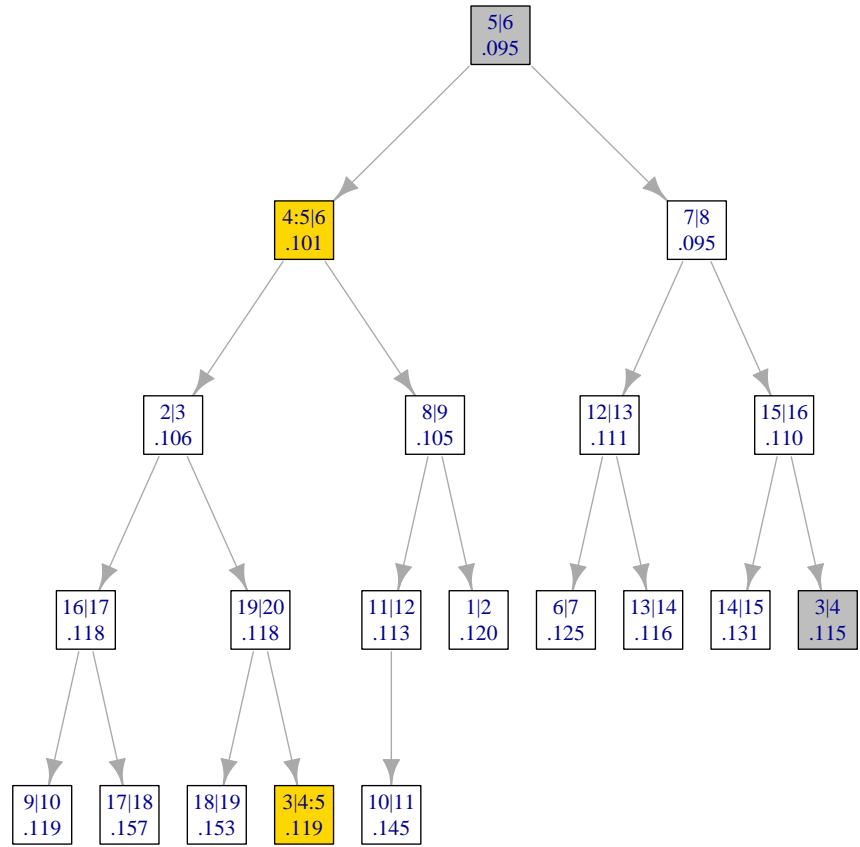

**Figure S3:** Min heap after the first merging step for the RLGH data set. The nodes corresponding to the fusion that have changed since initialization (Figure S1) are highlighted.

3. At each merging step the algorithm also tags as inactive the fusions involving the merged clusters (13). Indeed, once a cluster is fused with its left neighbor it can no longer be fused with its right neighbor and vice-versa. These fusions are highlighted in pink in Figure S1 and in gray (once tagged) in Figure S3. In order to avoid invalid fusions, each candidate fusion has an active/inactive label (represented by the gray highlight in Figure S3), and the when retrieving the next best candidate fusion (line 5), the min heap is first cleaned by deleting its root as long as it corresponds to an inactive fusion. In the course of the whole algorithm this additional cleaning step will at worst delete $2p$ roots in $\mathcal{O}(p \log(p))$.

4. The insertion instructions in Algorithm S2 indicate that the heap not only contains the value of the candidate fusions, but also the left and right clusters of each fusion, and the preceding and next candidate fusions in the order of the original objects to be clustered. In practice this side information is not actually stored in the heap, but in a dedicated array, together with the values of the corresponding linkage and the validity statuses of each candidate fusion. The heap only stores the index of each fusion in that array. The state of this array before and after the first fusion for the RLGH data set are given in Tables S1 and S2.

| left | right | prev | next | linkage | valid |
|---|---|---|---|---|---|
| 1 | 2 | NA | 2 | 0.121 | 1 |
| 2 | 3 | 1 | 3 | 0.106 | 1 |
| 3 | 4 | 2 | 4 | 0.115 | 1 |
| 4 | 5 | 3 | 5 | 0.095 | 1 |
| 5 | 6 | 4 | 6 | 0.095 | 1 |
| ⋮ | ⋮ | | | | ⋮ |
| 18 | 19 | 17 | 19 | 0.153 | 1 |
| 19 | 20 | 18 | NA | 0.118 | 1 |

**Table S1:** State of the array after initialization of the clustering for the RLGH data set, as in Figure S1.

| label | left | right | prev | next. | linkage | valid |
|---|---|---|---|---|---|---|
| 1—2 | 1 | 2 | NA | 2 | 0.121 | 1 |
| 2—3 | 2 | 3 | 1 | 20 | 0.106 | 1 |
| ~~3—4~~ | ~~3~~ | ~~4~~ | ~~2~~ | ~~4~~ | ~~0.115~~ | ~~0~~ |
| ~~4—5~~ | ~~4~~ | ~~5~~ | ~~3~~ | ~~5~~ | ~~0.095~~ | ~~0~~ |
| ~~5—6~~ | ~~5~~ | ~~6~~ | ~~4~~ | ~~6~~ | ~~0.095~~ | ~~0~~ |
| 6—7 | 6 | 7 | 21 | 7 | 0.125 | 1 |
| 7—8 | 7 | 8 | 6 | 8 | 0.096 | 1 |
| ⋮ | | | ⋮ | | | ⋮ |
| 18—19 | 18 | 19 | 17 | 19 | 0.153 | 1 |
| 19—20 | 19 | 20 | 18 | NA | 0.118 | 1 |
| 3—4:5 | 3 | 4:5 | 2 | 21 | 0.120 | 1 |
| 4:5—6 | 4:5 | 6 | 20 | 6 | 0.101 | 1 |

**Table S2:** State of the array after the first merge in the clustering for the RLGH data set, as in Figure S3.

# S3 Supplementary results

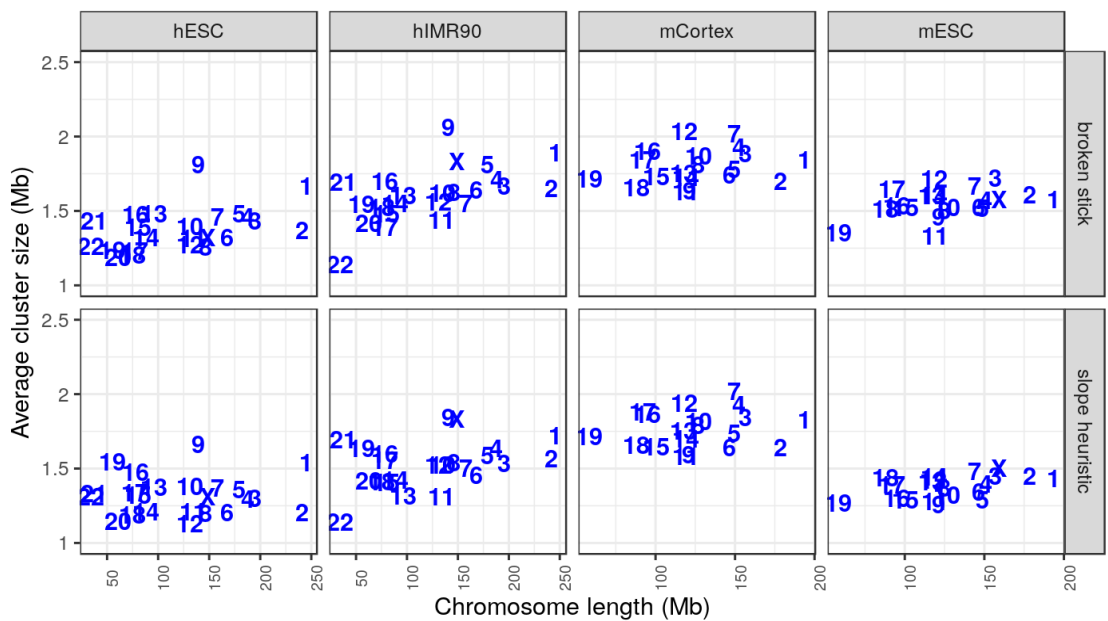

**Figure S4:** Average cluster size for both model selection approaches, compared to the chromosome length (in term of number of observed bins) for every chromosome and every experiment (full version). Chromosome X in mCortex had an average cluster size larger than 2.5Mb and was thus excluded from the picture.

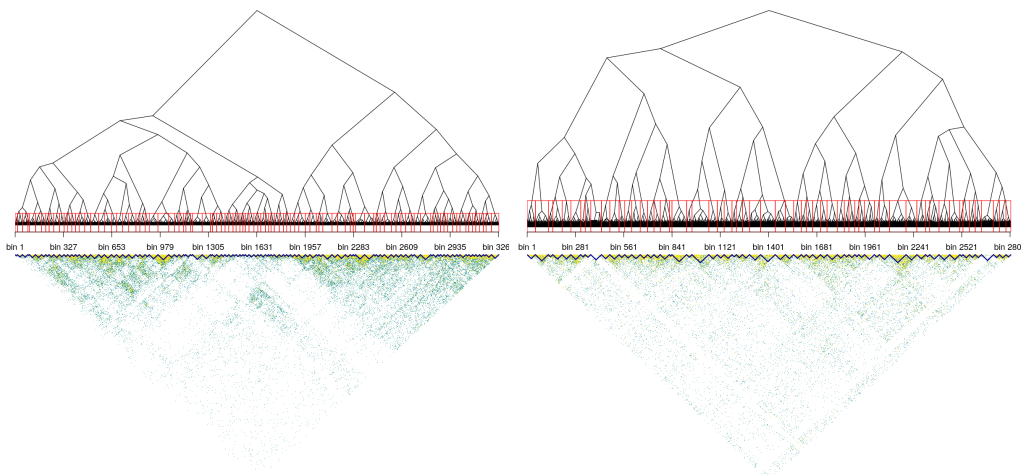

**Figure S5:** Left: Chromosome 11 of hIMR90. Right: Chromosome 12 of mCortex. Bottom: Hi-C data (log-scaled) with clustering selected by the slope heuristic (blue line). Top: Constrained hierarchical clustering with clustering selected by the slope heuristic (red rectangles).

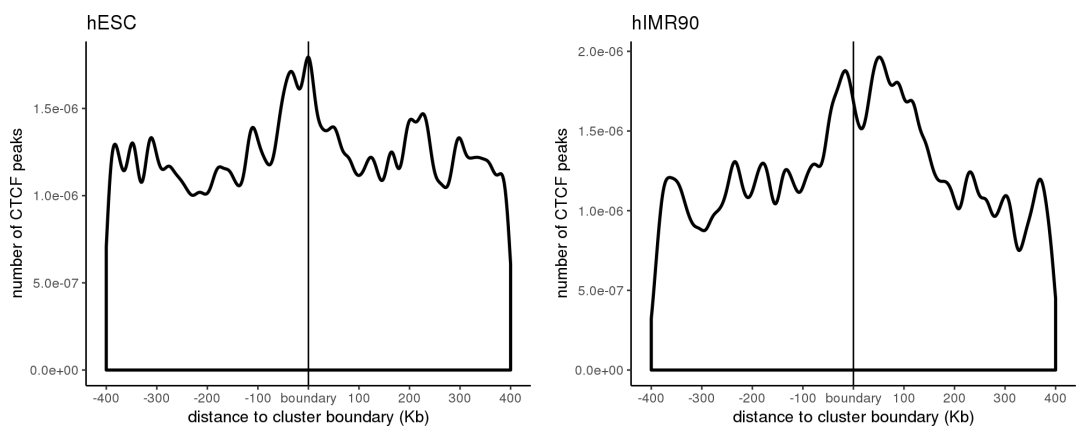

**Figure S6:** Distribution of the number of the 20% most intense CTCF ChIP-seq peaks with respect to distance of cluster boundaries, as obtained with the broken stick heuristic. Left: hESC. Right: hIMR90.

# References

J. Ah-Pine and X. Wang. Similarity based hierarchical clustering with an application to text collections. In H. Boström, A. Knobbe, C. Soares, and P. Papapetrou, editors, *Proceedings of the 15th International Symposium on Intelligent Data Analysis (IDA 2016)*, Lecture Notes in Computer Sciences, pages 320–331, Stockholm, Sweden, 2016. doi: 10.1007/978-3-319-46349-0. URL https://hal.archives-ouvertes.fr/hal-01437124.

N. Aronszajn. Theory of reproducing kernels. *Transactions of the American Mathematical Society*, 68(3):337–404, 1950.

9

F. Murtagh and P. Legendre. Ward's hierarchical agglomerative clustering method: which algorithms implement ward's criterion. *Journal of Classification*, 31:274–295, 2014. doi: 10.1007/s00357-014-9161-z.

J. Qin, D. P. Lewis, and W. S. Noble. Kernel hierarchical gene clustering from microarray expression data. *Bioinformatics*, 19(16):2097–2104, 2003. doi: 10.1093/bioinformatics/btg288.



\end{document}